\documentclass[12pt]{article}
\newtheorem{theorem}{Theorem}

\newtheorem{remark}[theorem]{Remark}

\newenvironment{proof}[1][Proof]{\textbf{#1.} }{\ \rule{0.5em}{0.5em}}
\newdimen\dummy
\dummy=\oddsidemargin
\begin{document}
\date{}
\title{Topological fundamental groups can distinguish spaces with isomorphic
homotopy groups}
\author{Paul Fabel \\
%EndAName
Department of Mathematics \& Statistics\\
Mississippi State University}
\maketitle

\begin{abstract}
We exhibit path connected metric spaces and a map $f:X\rightarrow Y$ such
that $f$ induces an isomorphism on homotopy groups but such that, with
natural topologies, $X$ and $Y$ have nonhomeomorphic fundamental groups.
Consequently we can conclude $X$ and $Y$ have distinct homotopy types
despite the failure of the Whitehead theorem in this context.
\end{abstract}

\section{Introduction}

Given CW complexes $X$ and $Y,$ the Whitehead theorem (\cite{hatch}) asserts
that a map $f:X\rightarrow Y$ is a homotopy equivalence provided $f$ induces
an isomorphism on homotopy groups. However the result can fail in the
context of path connected metric spaces. For example the standard Warsaw
circle has trivial homotopy groups but fails to have the homotopy type of a
point. This note aims to show the \textit{topological fundamental group }can
help counterbalance the general failure of the Whitehead theorem.

For a general space $X$ work of Biss \cite{Biss} initiates the development
of a theory whose fundamental notion is the following. Endowed with the
quotient topology inherited from the path components of based loops in $X$,
the familiar based fundamental group $\pi _{1}(X,p)$ of a topological space $%
X$ becomes a \textit{topological group}. For example if $X$ is locally
contractible then loops in $X$ are homotopically invariant under small
perturbation, and consequently the fundamental group $\pi _{1}(X,p)$ has the
discrete topology. For spaces that are complicated both locally and
globally, the topology of $\pi _{1}(X,p)$ can be more interesting (\cite{fab}
\cite{fab2} \cite{fab5}). An important feature of the theory is that if $X$
and $Y$ have the same homotopy type then $\pi _{1}(X,p)$ and $\pi _{1}(Y,p)$
are isomorphic \textit{and} homeomorphic (Proposition 3.3 \cite{Biss}).

These facts motivate the question of whether the added topological structure
on $\pi _{1}(X,p)$ can ever succeed in distinguishing the homotopy type of
spaces $X$ and $Y$ in instances when the hypotheses of the Whitehead theorem
are satisfied.

In fact this is so and in this note we exhibit aspherical spaces $X$ and $Y$
such that inclusion $j:X\rightarrow Y$ induces an isomorphism on homotopy
groups. However $\pi _{1}(X)$ and $\pi _{1}(Y)$ fail to be homeomorphic, and
thus we can conclude that $X$ and $Y$ do not have the same homotopy type
despite the failure of the Whitehead theorem for this pair of examples.

The theory of topological fundamental groups is still in the early stages of
development (\cite{fab3} \cite{fab4} \cite{fab6}) and it is hoped this note
will be seen as promoting its utility and helping to motivate its continued
investigation. For example the space $Y$ constructed in this paper is not
locally path connected. This suggests the following.

\textbf{Question. }Suppose $Y$ is an aspherical (metric) Peano continuum and 
$X\subset Y$ is aspherical and path connected. Suppose inclusion $%
j:X\hookrightarrow Y$ induces an isomorphism $j^{\ast }:\pi
_{1}(X,p)\rightarrow \pi _{1}(Y,p).$ Must $j^{\ast }$ be a homeomorphism? If 
$j^{\ast }$ is a homeomorphism must $j$ be a homotopy equivalence?

\section{Definitions and Preliminaries}

All definitions are compatible with those found in Munkres \cite{Munk}. If $%
X $ is a metrizable space and $p\in X$ let $C_{p}(X)=\{f:[0,1]\rightarrow X$
such that $f$ is continuous and $f(0)=f(1)=p\}.$ Endow $C_{p}(X)$ with the
topology of uniform convergence.

The \textbf{topological fundamental group} $\pi _{1}(X,p)$ is the set of
path components of $C_{p}(X)$ endowed with the quotient topology under the
canonical surjection $q:C_{p}(X)\rightarrow \pi _{1}(X,p)$ satisfying $%
q(f)=q(g)$ if and only if $f$ and $g$ belong to the same path component of $%
C_{p}(X).$ Thus a set $U\subset \pi _{1}(X)$ is open if and only if $%
q^{-1}(U)$ is open in $C_{p}(Y).$

\begin{remark}
The topological fundamental group $\pi _{1}(X,p)$ is a topological group
under concatenation of paths. (Proposition 3.1\cite{Biss}). A map $%
f:X\rightarrow Y$ determines a continuous homomorphism $f^{\ast }:\pi
_{1}(X,p)\rightarrow \pi _{1}(Y,f(p))$ via $f^{\ast }([\alpha ])=[f(\alpha
)] $ (Proposition 3.3 \cite{Biss}). If $X$ and $Y$ have the same homotopy
type then $\pi _{1}(X)$ is homeomorphic and isomorphic to $\pi _{1}(Y)$
(Corollary 3.4 \cite{Biss}) For the remainder of this paper all fundamental
groups will be considered topological groups.
\end{remark}

The space $X$ is \textbf{semilocally simply connected} at $p$ if there
exists an open set $U\subset X$ such that inclusion $j:U\hookrightarrow X$
induces the trivial homomorphism $j^{\ast }:\pi _{1}(U,p)\rightarrow \pi
_{1}(X,p).$ The space $Z$ is \textbf{discrete} if each one point subset of $%
Z $ is open.

\begin{remark}
\label{rem2}The main result of \cite{fab} shows that if $X$ is locally path
connected then $\pi _{1}(X,p)$ is discrete if and only if $\pi _{1}(X,p)$ is
semilocally simply connected.
\end{remark}

\section{Main result}

\begin{theorem}
There exist path connected aspherical separable metric spaces $X$ and $Y$
such that $X\subset Y$ and inclusion $j:X\hookrightarrow Y$ induces an
isomorphism $j^{\ast }:\pi _{1}(X,p)\rightarrow \pi _{1}(Y,p)$. Thus $%
(X,Y,j) $ satisfies the hypothesis of the Whitehead theorem. However the
topological fundamental groups $\pi _{1}(X,p)$ and $\pi _{1}(Y,p)$ are not
homeomorphic. Hence the \textbf{topology} of fundamental groups has the
capacity to distinguish the homotopy type of $X$ and $Y$ when the algebra
fails to do so.
\end{theorem}

\begin{proof}
The basic idea is to let $X$ denote the countable union of a sequence of
large simple closed curves $C_{1}\cup C_{2}...$ joined at a common point $p.$
Such a space is sometimes called a bouquet of infinitely many loops. In
particular $X$ is locally contractible and should not be mistaken for the
Hawaiian earring. The space $Y$ is a compactification of $X$ obtained by
attaching a line segment $\alpha $ based at $p$ such that the curves $C_{n}$
converge to $\alpha $ in the Hausdorff metric.

Since each of $X$ and $Y$ is path connected and 1 dimensional, if $n\neq 1$
then $\pi _{n}(X,p)=\pi _{n}(Y,p)=1.$ Thus, to show that $(X,Y,j)$ satisfies
the hypothesis of the Whitehead theorem it suffices to show that $j^{\ast
}:\pi _{1}(X,p)\rightarrow \pi _{1}(Y,p)$ is an isomorphism.

Formally for $n\geq 2$ let $C_{n}\subset R^{2}$ denote boundary of the
convex hull of the following 3 point set: $\{(0,0),(\frac{1}{n},1),(\frac{1}{%
n},1)+\frac{1}{10^{(10n)}}(n,-1)\}.$ Then for each $n\geq 2$ $C_{n}$ is the
boundary of a triangle and in particular $C_{n}$ is a simple closed curve.
Let $p=(0,0).$ Note $C_{n}\cap C_{m}=p$ if $n\neq m.$ Let $\alpha $ denote
the line segment $[(0,0),(0,1)]\subset R^{2}.$ Let $X=\cup _{n=2}^{\infty
}C_{n}$ and let $Y=\overline{X}.$ Note $X\cup \alpha .$ Note the path
connected spaces $X$ and $Y$ are 1 dimensional and hence aspherical (\cite
{Fort}). We will show inclusion $j:X\hookrightarrow Y$ induces an
isomorphism $j^{\ast }:\pi _{1}(X,p)\rightarrow \pi _{1}(Y,p).$

To prove $j^{\ast }$ is one to one suppose $f:\partial D^{2}\rightarrow X$
is inessential in $Y$ and suppose $f(1)=p.$ Let $F:D^{2}\rightarrow Y$
satisfy $F_{\partial D^{2}}=f.$ Let $U=F^{-1}(\alpha \backslash p).$ Since $%
D^{2}$ is locally path connected, and since $\alpha \backslash p$ is a
component of $X\backslash p$ the set $U$ is open. Suppose $x\in \overline{U}%
\backslash U.$ Then $F(x)=p$ since $\alpha =\overline{\alpha \backslash p}.$
Thus, we may redefine $F$ to be $p$ on the set $U$ and obtain a continuous
function $G:D^{2}\rightarrow X$ such that $G_{\partial D^{2}}=f.$ This
proves $j^{\ast }$ is one to one.

To prove $j^{\ast }$ is a surjection suppose $\beta \in C_{p}(Y).$ We must
show there exists $\gamma \in C_{p}(X)$ such that $\gamma $ and $\beta $ are
path homotopic in $Y.$ Since $im(\beta )$ is a Peano continuum $im(\beta )$
is locally path connected. Thus we may choose $N$ such that $im(\beta )\cap
(\{\frac{1}{N},\frac{1}{N+1},..\}\times \{1\})=\emptyset .$ Let $A=im(\beta
)\cap (\alpha \cup C_{N}\cup C_{N+1}...).$ Let $B=C_{1}\cup C_{2}...\cup
C_{N-1}.$ Note $A$ is a contractible Peano continuum such that $p\in A.$
Moreover $B$ is a strong deformation retract of $B\cup A.$ Thus there exists
a homotopy $h_{t}:A\cup B\rightarrow B$ such that $h_{0}=id_{A\cup B}$ and $%
h_{t}$ fixes $B$ pointwise. Thus the homotopy $h_{t}(\beta )$ determines
that $\beta $ is path homotopic in $Y$ to $\gamma =h_{1}(\beta ).$ Note $%
im(\gamma )\subset X.$ Hence $j^{\ast }$ is a surjection and therefore $%
j^{\ast }$ is an isomorphism.

Since the space $X$ is locally contractible $\pi _{1}(X,p)$ has the discrete
topology (Remark \ref{rem2}). On the other hand $\pi _{1}(Y,p)$ does not
have the discrete topology, since there exists an inessential loop $f\in
C_{p}(Y)$ which is the uniform limit of inessential loops. (Let the
(inessential) map $f$ go up and down once on $\alpha $ and let $f_{n}$ be an
(essential) loop going once around $C_{n}$). Thus the path component of the
constant map is not open in $C_{p}(Y)$ and thus $\pi _{1}(Y,p)$ cannot have
the discrete topology. Thus $\pi _{1}(X,p)$ and $\pi _{1}(Y,p)$ are not
homeomorphic and hence $X$ and $Y$ do not have the same homotopy type.
\end{proof}

\begin{remark}
The space $Y$ constructed is semilocally simply connected. However $\pi
_{1}(Y,p)$ does not have the discrete topology. Consequently $Y$ is a
counterexample to the (false) Theorem 5.1 \cite{Biss} which asserts that $%
\pi _{1}(Y,p)$ is discrete if and only if $\pi _{1}(Y,p)$ is semilocally
simply connected.
\end{remark}


\begin{thebibliography}{99}
\bibitem{Biss}  Biss, Daniel K. \textit{The topological fundamental group
and generalized covering spaces.} Topology Appl. 124 (2002), no. 3, 355--371.

\bibitem{Fort}  Curtis, M. L.; Fort, M. K., Jr. \textit{Homotopy groups of
one-dimensional spaces.} Proc. Amer. Math. Soc. 8 (1957), 577--579.

\bibitem{fab}  Fabel, Paul. \textit{A characterization of spaces with
discrete topological fundamental group.} Preprint.
http://front.math.ucdavis.edu/math.GN/0502249

\bibitem{fab2}  Fabel, Paul. \textit{The fundamental group of the harmonic
archipelago.} Preprint.http://front.math.ucdavis.edu/math.AT/0501426

\bibitem{fab3}  Fabel, Paul. \textit{The topological Hawaiian earring group
does not embed in the inverse limit of free groups.} Preprint.
http://front.math.ucdavis.edu/math.GN/0501482

\bibitem{fab4}  Fabel, Paul \textit{A retraction theorem for topological
fundamental groups with applications to the Hawaiian earring. }%
Preprint.http://front.math.ucdavis.edu/math.AT/0502218

\bibitem{fab5}  Fabel, Paul \textit{The Hawaiian earring group is
topologically incomplete. }Preprint.
http://front.math.ucdavis.edu/math.GN/0502148

\bibitem{fab6}  Fabel, Paul \textit{A monomorphism theorem for the inverse
limit of nested retracts. }Preprint.
http://front.math.ucdavis.edu/math.AT/0502275

\bibitem{hatch}  Hatcher, Allen. \textit{Algebraic topology}. Cambridge
University Press, Cambridge, 2002.

\bibitem{Munk}  Munkres, James R., \textit{Topology: a first course.}
Prentice-Hall, Inc., Englewood Cliffs, N.J., 1975.
\end{thebibliography}
\end{document}